\documentclass{ifacconf}

\usepackage{mathtools}
\usepackage{graphicx}      
\usepackage{natbib}        
\usepackage{amsmath,amssymb,amsfonts}
\usepackage{mathrsfs}
\usepackage{array}
\usepackage{graphicx}
\usepackage{textcomp}
\usepackage{xcolor}
 \usepackage{enumerate}

\usepackage{algorithm}
\usepackage{algorithmic}
\newtheorem{assumption}{Assumption}

\newtheorem{remark}{Remark}

\newcommand{\ADD}[1]{\textcolor{black}{{#1}}}

\begin{document}
	\begin{frontmatter}
		
        \title{Lightweight Real-Time ALADIN for Distributed Optimization}
		

\thanks{$^\dagger$ Corresponding author}
\thanks[footnoteinfo]{The work of X.D. and A.I.R. was supported by the Guangzhou-HKUST(GZ) Joint Funding Scheme (Grant No. 2025A03J3960). The work of A.I.R. was also supported by the Guangdong Provincial Project (Grant No. 2024QN11G109).}
        
\author[First]{Yifei Wang}
\author[Second]{Xuhui Feng}
\author[Second]{Shimin Pan}
\author[Second]{Liangfan Zhu}
\author[Third]{Xu Du$^\dagger$} 
\author[Third,Forth]{Apostolos I. Rikos}
 
\address[First]{The Ningbo Artificial Intelligence Institute and the Department of Automation, Shanghai Jiao Tong University, Ningbo, China\\The School of Automation and Intelligent Sensing, Shanghai Jiao Tong University, Shanghai, China (e-mail: yifeiw4ng@sjtu.edu.cn)} 
\address[Second]{Huawei Technologies Co., Ltd., Shanghai, China (e-mail: \{fengxuhui, panshimin1, zhuliangfan\}@huawei.com)}
\address[Third]{The Artificial Intelligence Thrust of the Information Hub, The Hong Kong University of Science and Technology (Guangzhou), Guangzhou, China (e-mail: \{michaelxudu,apostolosr\}@hkust-gz.edu.cn).}
\address[Forth]{The Department of Computer Science and Engineering, The Hong Kong University of Science and Technology, Clear Water Bay, Hong Kong, China.}

\begin{abstract}
This paper presents a real-time computational framework for multi-node distributed optimization by extending the Augmented Lagrangian Alternating Direction Inexact Newton (ALADIN) algorithm. Our approach integrates adjoint sequential quadratic programming (SQP) techniques to enable efficient approximation of Jacobian information within the ALADIN embedded quadratic program, thereby reducing communication overhead. Furthermore, to decrease computational complexity, we design an event-triggered update strategy that avoids updating Hessian and Jacobian matrices at every iteration. 
The proposed method achieves local convergence and enhanced communication efficiency, making it well suited for time-critical applications.
Numerical experiments demonstrate that our approach achieves competitive performance while exhibiting superior computational efficiency in real-time scenarios, validating its practical applicability for time-sensitive distributed optimization challenges.
\end{abstract}
		
\begin{keyword}
Distributed optimization,  ALADIN, Adjoint SQP
\end{keyword}
		
\end{frontmatter}
	
\section{Introduction}\label{sec: intro}

{\ADD{Distributed optimization has attracted substantial attention in recent years. This trend is driven by progress in federated learning (\cite{ren2025communication,Du2023_Arxiv}), power systems (\cite{Du2019,lanza2025distributed}), and model predictive control (MPC)  (\cite{stomberg2025decentralized,Kouzoupis2016AA}). Studies in these fields consistently show that the curse of dimensionality becomes severe as the dimension of the optimization variables grows. As a result, centralized methods for large-scale problems become computationally infeasible. Distributed optimization algorithms are therefore viewed as an effective means to alleviate the curse of dimensionality.}}

{\ADD{Distributed optimization comprises a family of computational methods in which nodes iteratively solve local subproblems. During the iteration process, nodes exchange information with neighbors or with a central coordinator. The exchanged information may include first and second-order information, or local primal variables. These approaches are commonly classified into primal and dual decomposition frameworks (\cite{2024_doostmohammadian_rikos_Johansson_survey}). In dual decomposition methods, communication includes not only the information exchanged in primal decomposition but also dual variables associated with coupling constraints, which can enhance convergence behavior.
A well-known algorithm in this category is the Alternating Direction Method of Multipliers (ADMM) (\cite{boyd2011distributed}), which distributes the Augmented Lagrangian Method and provides rigorous convergence guarantees for convex problems. However, ADMM often suffers from slow convergence and requires restrictive assumptions in nonconvex settings, see \cite{hong2016convergence}. Recent research has introduced the Augmented Lagrangian Alternating Direction Inexact Newton algorithm (ALADIN) (\cite{Shi2022,Du2025}). ALADIN incorporates sequential quadratic programming (SQP) into an ADMM-type framework, which yields accelerated convergence for convex problems and ensures local convergence for nonconvex cases. Building on these properties, this paper aims to explore the application of ALADIN to distributed nonconvex optimization in real-time environments.}}


\noindent
\textbf{Existing Literature.}
 {\ADD{ALADIN was originally developed for distributed sensor network localization (\cite{Houska2016}) and MPC (\cite{Kouzoupis2016AA}), both of which involve distributed nonconvex optimization problems. Subsequent studies (\cite{Engelmann2019,Du2019}) demonstrated its effectiveness in distributed alternating current power flow optimization. All of these applications fall within the category of resource allocation problems. More recently, \cite{Du2025} introduced a consensus-oriented variant, termed Consensus ALADIN, which achieves improved computational efficiency for consensus optimization. To the best of our knowledge, the paper (\cite{wangCDC2025}) was the first to extend  the application of ALADIN to non-smooth, non-convex constraint optimization problems.}}
Regarding convergence properties, ALADIN ensures local convergence for nonconvex problems under standard assumptions. These assumptions include the second order sufficient condition (SOSC) and the linear independent constraint qualification (LICQ). For convex problems, global convergence can be achieved when the local Hessian is approximated by constant matrices. Recent theoretical advances (\cite{Du2023_Arxiv, Du2025, Du2025CDCA}) established sufficient conditions for global linear convergence of ALADIN in convex settings. Moreover, \cite[Chapter~5]{ALADINtheisi} derived explicit linear convergence rates under stronger assumptions. For comprehensive surveys of ALADIN, we refer interested readers to \cite{lanza2025distributed}.

However, a key limitation of the ALADIN framework lies in its reliance on solving a large-scale constrained coupled quadratic program (QP) for inter-node coordination. This step introduces high communication overhead and significant computational cost, mainly due to the construction and direct solution of the coupled QP (\cite[IV.C]{Du2019}). Recent work (\cite{Wu2025CDC}) derived a closed-form expression for the constrained coupled QP. This result decouples the updates of three variable sets: (i) the dual variables associated with coupling constraints, (ii) the dual variables related to local constraints, (iii) the primal variables.
Despite this simplification, constructing the coordination QP still requires each node to transmit substantial local information to the coordinator when local constraints are present. This information includes first and second-order information of the local objectives, the constraint Jacobians, which results in significant communication overhead. Moreover, the variable updates in the closed-form representation still rely on matrix inversions. These inversions become computationally expensive in high-dimensional settings.
Although existing Hessian approximation techniques, such as the Broyden–Fletcher–Goldfarb–Shanno (BFGS)-based method in \cite{Du2025}, can avoid the direct transmission of Hessian information, they are only applicable to unconstrained subproblems and do not provide a general mechanism for handling constrained cases. These limitations reveal two major challenges for ALADIN:
(i) excessive information transmission and the absence of effective Hessian and Jacobian recovery for constrained local subproblems, and
(ii) considerable matrix inversion operations required to solve the large-scale coupled QP.

\noindent
\textbf{Main Contributions.} Motivated by the aforementioned challenges, we propose a novel real-time optimization framework based on ALADIN, termed \emph{Efficient Adjoint BFGS ALADIN}. Our approach combines the BFGS-ALADIN technique (\cite{Du2025}) with adjoint SQP methods 
(\cite{diehl2010adjoint,quirynen2018inexact,feng2020inexact,hespanhol2021adjoint}),
enabling the coordinator to reconstruct both local Hessian and Jacobian information without explicit transmission from the nodes. Furthermore, inspired by \cite{Shi2022}, we incorporate an event-triggered mechanism that adaptively determines when to update these local approximations during the iterative process. To the best of our knowledge, this is the first work that employs ALADIN within a framework that simultaneously (i) recovers both local Hessian and Jacobian information at the coordinator and (ii) integrates an event-triggered mechanism to update these approximations, thereby reducing the computational burden of solving ALADIN's coupled QP. Numerical results demonstrate the superior performance of the proposed algorithm compared to existing benchmarks.

\emph{Notation.} In this paper, the superscript $[k]$ denotes the iteration index, while the notation $(\psi\mid \phi)$ indicates the dual multiplier $\psi$ associated with the corresponding constraint $\phi$.
    \section{Problem Formulation} \label{sec: OPT problem}
This paper considers a resource allocation problem involving \(N\) nodes. Each node \(i\) has a decision variable \(x_i \in \mathbb{R}^{n_i}\), a local objective function \(f_i(x_i): \mathbb{R}^{n_i} \rightarrow \mathbb{R}\) and local constraints \(g_i: \mathbb{R}^{n_i} \rightarrow \mathbb{R}^{c_i}\). The problem is formulated as follows:
\begin{equation}\label{eq: DIS OPT}\small
    \begin{split}
        \min_{\substack{x_i\in \mathbb R^{n_i} \\ \forall i \in \{1,\dots,N\}}} \; \sum_{i=1}^N f_i(x_i)\quad  \text{s.t.}\;\;\;\left\{\begin{split}
             &  \sum_{i=1}^N  A_ix_i = b, \\
        &  g_i(x_i) = 0\quad \forall i \in \{1,\dots,N\}.
        \end{split}\right.
    \end{split}
\end{equation}
Here, the nodes' variables are coupled through an affine constraint \(\sum_{i=1}^N A_i x_i = b\), where \(A_i \in \mathbb{R}^{m \times n_i}\) and \(b \in \mathbb{R}^m\). 
   
In this paper, we adopt the following assumptions that are essential for our subsequent development.
\begin{assumption}\textbf{\emph{SOSC Satisfaction}}.\label{ass 1}
    For each node, the local cost function $f_i$ is closed, proper, twice continuously differentiable, and satisfies SOSC at a local minimizer $x^*$ of problem \eqref{eq: DIS OPT}, i.e. $\nabla^2 f(x^*) \succ 0.$
    Additionally, for every node $i$, the local constraint function $g_i$ is assumed to be closed, proper, and second-order continuously differentiable.
\end{assumption}

\begin{assumption}\textbf{\emph{LICQ Satisfaction}}.\label{ass 2}
At the local minimizer $x^\star$ of problem \eqref{eq: DIS OPT}, the Jacobian of all equality constraints is full row rank. This ensures that the KKT system associated with problem \eqref{eq: DIS OPT} is regular.
\end{assumption}

Assumptions \ref{ass 1} and \ref{ass 2} guarantee that $x^\star$ is a regular Karush-Kuhn-Tucker (KKT) point of problem \eqref{eq: DIS OPT}, which forms the basis for the local convergence of ALADIN-type methods for nonconvex problems. For further details, we refer to \cite{Houska2016,Kouzoupis2016AA}.

    \section{Preliminaries: Efficient ALADIN}\label{sec: ALADIN}


The Lagrangian function of problem \eqref{eq: DIS OPT} is the following,
\begin{equation}\label{eq: largangian}\small
    \mathcal L(x, \lambda) =\left( \sum_{i=1}^N f_i(x_i) + \mu_i^\top g_i(x_i) \right)+ \lambda^\top \left(\sum_{i=1}^N A_i x_i-b\right).
\end{equation}
Here, ${x} = [x_1^\top, x_2^\top, \cdots, x_N^\top]^\top$, where $\mu_i$ denotes the local dual variable associated with the constraint $g_i(x_i)=0$, and $\lambda$ indicates the coupled dual variable corresponding to the coupling constraint. 

Based on \eqref{eq: largangian}, 
in our previous work \cite{Wu2025CDC}), we developed an \emph{Efficient Gauss-Newton ALADIN} framework for solving moving horizon estimation (MHE) problems. The transformation from the distributed formulation in \eqref{eq: DIS OPT} to the algorithmic structure in \eqref{eq: EGN-ALADIN} follows this framework, and the detailed derivation of the augmented Lagrangian decomposition and the associated coordination QP can be found in \cite{Wu2025CDC}. This approach reformulates the centralized MHE problem as a distributed optimization problem in the form of \eqref{eq: DIS OPT} (details can be found in Section \ref{subsection: distributed model}). The algorithmic structure is presented below:
\begin{equation}\label{eq: EGN-ALADIN}\small
\left\{
    \begin{split}
     &(x_i^{[k+1]}, \mu_i^{[k+1]})= \arg\min_{x_i}  f_i(x_i) + (\lambda^{[k]})^\top x_i + \frac{\rho}{2} \|x_i - y_i^{[k]}\|^2,  \\
        & \qquad \qquad\qquad\text{s.t.} \quad  g_i(x_i) = 0\;|\;\mu_i,\quad \forall i \in \mathcal{V},\\
        &\; H_i \approx \nabla^2 (f_i(x_i^{[k+1]})+ ( \mu_i^{[k+1]})^\top g_i(x_i^{[k+1]})) \succ 0,\\
        &\; v_i = \nabla f_i(x_i^{[k+1]}),\; C_i = \nabla g_i(x_i^{[k+1]}), \\
           & \;\lambda^{[k+1]} = \left(\mathop{\sum}_{i=1}^{N}G_i - Q_iR_i^{-1}Q_i^\top \right)^{-1}p,\\
               &\;\tilde{\mu_i}=  -R_i^{-1}\left(C_iH_i^{-1}v_i+Q_i^\top \lambda^{[k+1]}\right),\\
	&\;\Delta x_i=  -H_i^{-1}\left(\nabla f_i(x_i^{[k+1]})+C_i^\top \tilde{\mu_i}+A_i^\top \lambda^{[k+1]}\right),
       \\
        &\; y^{[k+1]} = x^{[k+1]} + \Delta x.
    \end{split}
\right.
\end{equation}
\vspace{-1mm}
where the auxiliary matrices and vectors are defined as:
\begin{equation}\label{eq: parameters}\small
			\left\{
				\begin{split}
                \hspace{-1mm}G_i\hspace{-0.8mm}=&A_iH_i^{-1}A_i^\top, \\
                   \hspace{-1mm}Q_i\hspace{-0.8mm}=&A_iH_i^{-1}C_i^\top, \\
                   \hspace{-1mm}R_i\hspace{-0.8mm}=&C_iH_i^{-1}C_i^\top,
                \end{split}
            \right.
            \quad
            \text{and}
            \quad
            \left\{
            	\begin{split}
                q\hspace{-0.8mm}=\hspace{-0.8mm} &\mathop{\sum}_{i=1}^{N}\left(Q_iR_i^{-1}C_i-A_i\right)H_i^{-1}v_i,\\
                   p\hspace{-0.8mm}=\hspace{-0.8mm}&\mathop{\sum}_{i=1}^{N}A_ix_i^{[k+1]}+q.
                \end{split}
			\right.
\end{equation}
Moreover, let ${y} = [y_1^\top, y_2^\top, \cdots, y_N^\top]^\top$ denote the stacked global variables and $\Delta x = [\Delta x_1^\top, \Delta x_2^\top, \cdots, \Delta x_N^\top]^\top$ their corresponding increments {\ADD{(see \eqref{eq: EGN-ALADIN}).}} For each subproblem, $H_i$, $C_i$, and $v_i$ denote the Gauss-Newton Hessian approximation, the constraint Jacobian, and the local objective gradient, respectively. The variables $\tilde{\mu}_i$ represent the dual multipliers of the linearized local constraints in the coordination QP and are used only to compute the primal search directions, unlike the local multipliers $\mu_i^{[k+1]}$ obtained from the local subproblems.

The algorithm consists of three main steps: first, each node solves its local augmented Lagrangian subproblem with parameter $\rho$ and the previously updated primal variable $y_i^{[k]}$ received from all nodes; second, based on the updated local solution $x_i^{[k+1]}$, the gradient $\upsilon_i=\nabla f_i(x_i^{[k+1]})$, Hessian $H_i$ of each objective function $f_i$, and the local constraint Jacobian $C_i$ of $g_i(x_i^{[k+1]})$ are evaluated; third, the coordinator updates $\lambda^{[k+1]}$, $\tilde
\mu_i$, and $\Delta x_i$ for all $i \in \{1,2,\cdots, N\}$ using the closed-form expression for the affine-coupled constrained QP proposed in \cite[Appendix I]{Wu2025CDC}):
\begin{equation}\label{eq: coupled QP}\small
    \begin{split}
        \min_{\substack{\Delta x_i\in \mathbb R^{n_i} \\ \forall i \in \{1,\dots,N\}}} &\quad \sum_{i=1}^N \frac{1}{2}\Delta x_i^\top H_i \Delta x_i + \upsilon_i^\top \Delta x_i \\
        \text{s.t.} \quad& \quad \sum_{i=1}^N A_i(\Delta x_i + x_i^{[k+1]}) = b \;|\; \lambda \\
        & \qquad C_i\Delta x_i = 0, \;|\; \tilde \mu_i \quad \forall i \in \{1,\dots,N\},
    \end{split}
\end{equation}
and broadcast $(\lambda^{[k+1]}, y^{[k+1]})$ to all nodes, see the last step of \eqref{eq: EGN-ALADIN}. Repeat the above steps until convergence.
Note that Assumptions \ref{ass 1} and \ref{ass 2} guarantee the existence and uniqueness of the solution to the coupled QP \eqref{eq: coupled QP}.

\section{LIGHTWEIGHT ALADIN}
We now present the main result of this paper.
Specifically, in this section we address two major challenges of ALADIN in resource allocation problems: (i) the high communication load caused by transmitting Hessian and Jacobian information, and (ii) the computational cost of solving large-scale coupled QPs. To address these issues, we integrate an adjoint-based BFGS scheme that allows the coordinator to reconstruct accurate Hessian and Jacobian approximations without explicit transmission. Furthermore, we introduce an event-triggered update scheme, which keeps local matrices fixed between triggers. This approach significantly improves computational efficiency while preserving convergence.
  \subsection{Efficient Adjoint BFGS ALADIN}\label{sec: Adjoint ALADIN} 
\begin{algorithm}[ht]\small
	\caption{Efficient Adjoint BFGS ALADIN}
	\textbf{Initialization:} Initial guess of dual variable $\lambda$ and primal variables $\{Y_i\}, \;\forall i$, choose $\rho>0$. \\
    \textbf{Output:} Optimal solution $\{Y_i^\star\}$. \\
	\textbf{Repeat:}
	\begin{enumerate}
		\item[1.] Paralleled solve local NLP:
		\begin{equation}\label{ALADIN-step1} 
        \begin{aligned}
\hspace{-0.5mm}({x_i}^{[k+1]},\;\mu_i^{[k+1]})
\hspace{-0.5mm}=\hspace{-0.5mm}\mathop{\arg\min}_{x_i} \hspace{0.5mm}&f_i(x_i)\hspace{-0.5mm}+(\lambda^{[k]})^\top \hspace{-0.5mm}A_iX_i\hspace{-0.4mm}+\frac{\rho}{2}\|x_i\hspace{-0.5mm}-\hspace{-0.5mm}y_i^{[k]}\hspace{-0.4mm}\|^2\\
        \mathrm{s.t.} \quad & g_i(x_i)=0\;|\;\mu_i.
        \end{aligned}
		\end{equation}
        
		\item[2.]  \ADD{Evaluate the local gradient $\upsilon_i^{[k+1]} = \nabla f_i(x_i^{[k+1]})$, and update the auxiliary quantities required for the Hessian and Jacobian  recovering:
\begin{equation}\label{auxiliary vectors updates}\left\{
    \begin{aligned}
        S_i &= x_i^{[k+1]} - x_i^{[k]},\quad d_i = \upsilon_i^{[k+1]} - \upsilon_i^{[k]},\\
         z_i &= g_i(x_i^{[k+1]}) - g_i(x^{[k]}_i),\quad \sigma_i = \mu_i^{[k+1]} - \mu_i^{[k]}.
    \end{aligned}\right.
\end{equation}
Then compute the adjoint direction using the reverse mode}
\begin{equation}\label{gamma updates}
    \gamma_i^{[k+1]} = \nabla_{x_i}g_i(x_i^{[k+1]})^\top \sigma_i^{[k+1]}.
\end{equation}
\item[3.]
Update the BFGS Hessian approximation as
\begin{equation}\label{eq:BFGS hessian update}
    H_i^{[k+1]} = H_i^{[k]} -\frac{H_i^{[k]}S_iS_i^{\top} H_i^{[k]}}{S_i^{\top} H_i^{[k]} S_i} + \frac{d_i d_i^{\top}}{S_i^{\top} d_i}.
\end{equation}
Moreover, the Jacobian approximation updates as showed in \cite{diehl2010adjoint}
\begin{equation}\label{eq: adjoint}
    C_i^{[k+1]} = C_i^{[k]} +\frac{\left(z_i -C_i^{[k]} S_i \right)\left( (\gamma_i^{[k+1]})^\top -\sigma_i^{\top}C_i^{[k]} \right)}{\left( (\gamma_i^{[k+1]})^\top -\sigma_i^{\top}C_i^{[k]} \right) S_i}.
\end{equation}
		\item[4.] Update the global dual variable $\lambda$ with the parameters in \eqref{eq: parameters}:
		 \begin{equation}\label{eq:global dual variable}
                \lambda^{[k+1]} =\left(\mathop{\sum}_{i=1}^{N}G_i - Q_iR_i^{-1}Q_i^\top \right)^{-1}p.
		\end{equation}  
         Paralleled update local primal and dual variables:
		\begin{equation}\label{eq:local variables}
			\left\{
				\begin{aligned}
					\tilde{\mu_{i}}= & -R_i^{-1}(C_i^{[k+1]}H_i^{-1}\upsilon_i+Q_i^\top \lambda^{[k+1]}),\\
					y_{i}^{[k+1]}= &\;x_i^{[k+1]} -H_i^{-1}(\upsilon_i+(C_i^{[k+1]})^\top \tilde{\mu_{i}}+A_i^\top \lambda^{[k+1]}).
				\end{aligned}
			\right.
		\end{equation} 
	\end{enumerate}
	\label{alg:ALADIN}
\end{algorithm}
Although the above algorithm can solve the distributed optimization problem \eqref{eq: DIS OPT} with high accuracy, it requires substantial uplink and downlink communication. This section proposes an adjoint-based BFGS scheme to address this issue. The proposed algorithm, referred to as \emph{Efficient Adjoint BFGS ALADIN}, is summarized in Algorithm \ref{alg:ALADIN}.

Following the algorithmic framework in \eqref{eq: EGN-ALADIN}, the proposed method consists of four steps: (i) Following the first step of \eqref{eq: EGN-ALADIN}, each node solves its local NLP subproblem to compute the updated primal variable $x_i^{[k+1]}$ and the associated multiplier $\mu_i^{[k+1]}$. (ii) Based on $(x_i^{[k+1]}, \mu_i^{[k+1]})$, it computes the local gradient $\upsilon_i^{[k+1]}$ and updates auxiliary vectors $(S_i,d_i,z_i,\sigma_i,\gamma_i^{[k+1]})$ as defined in \eqref{auxiliary vectors updates} and \eqref{gamma updates}. The updated quantities $(S_i,d_i, z_i, \sigma_i, \gamma_i^{[k+1]})$ are then transmitted to the coordinator\footnote{
Transmitting either variable differences or the variables themselves incurs the same communication load. This work adopts the former approach, where variable differences are computed in parallel across all nodes.
}. (iii) Using the adjoint-based BFGS update described in \cite{diehl2010adjoint} and shown in \eqref{eq:BFGS hessian update} and \eqref{eq: adjoint}, the coordinator refines the Hessian and Jacobian approximations with high accuracy. (iv) Finally, the dual variables $\lambda^{[k+1]}$ and $\tilde{\mu_i}$ together with the primal direction $\Delta x_i$ are updated according to \eqref{eq:global dual variable} and \eqref{eq:local variables} for all $i \in \{1,2,\cdots, N\}$.


\begin{remark}
    \textbf{\emph{Reduction in Communication Overhead}}. The primary distinction from \eqref{eq: EGN-ALADIN} lies in the third step, where the proposed algorithm avoids both explicit construction and transmission of full sensitivity matrices (Hessian and Jacobian), thereby reducing computational and communication overhead. In \eqref{eq: EGN-ALADIN}, each node transmits the tuple $(x_i^{[k+1]}, \upsilon_i^{[k+1]}, H_i, C_i)$ to the coordinator, yielding a communication load of $(2n_i + n_i^2 + c_i \times n_i)$. For least-squares problems under the \emph{Gauss-Newton ALADIN} scheme, the Hessian transmission is replaced by that of the measurement Jacobian, as discussed in \cite{Du2019}. In contrast, Algorithm \ref{alg:ALADIN} reconstructs accurate Hessian and Jacobian approximations at the coordinator using the adjoint-based BFGS update. As a result, each node only transmits  $(S_i,d_i,z_i,\sigma_i,\gamma_i^{[k+1]})$, as defined in \eqref{auxiliary vectors updates} and \eqref{gamma updates}, reducing the communication load to $(2n_i+3c_i)$.
\end{remark}

As shown in \eqref{eq: parameters}, the coordinator update involves repeated matrix inversions, which remain the main computational bottleneck. This issue is alleviated in the following subsection.

\subsection{Real-Time Efficient Adjoint BFGS ALADIN}\label{sec: Real Time Adjoint ALADIN} 
This subsection introduces an event-triggered update scheme to accelerate Algorithm \ref{alg:ALADIN}. Inspired by \cite[Algorithm 1, Step 9]{Shi2022}, the Hessian and Jacobian approximations $H_i^{[k]}$ and $C_i^{[k]}$, together
with the associated matrices $G_i$, $Q_i$, and $R_i$ defined in \eqref{eq: parameters}, are updated only when $\log_3(k) \in \mathbb{N}$ and kept constant otherwise. 
This design is motivated by the observation that, as the iterates $x_i^{[k+1]}$ and $x_i^{[k]}$ gradually converge, the sensitivity matrices gradually stabilize, making frequent updates unnecessary. Between two triggering iterations, only the local NLP solves in \eqref{ALADIN-step1} and the coordination updates in \eqref{eq:global dual variable}–\eqref{eq:local variables} are performed, while the sensitivity matrices in \eqref{eq:BFGS hessian update}–\eqref{eq: adjoint} are reused. Hence, the proposed heuristic reduces the update frequency over time without compromising convergence performance. Unlike \cite{Shi2022}, we extend this event-triggered sensitivity strategy to the nonconvex setting considered in this paper. Moreover, by exploiting the closed-form solution of the coordination QP in \eqref{eq:global dual variable}–\eqref{eq:local variables}, the proposed approach eliminates the need for external QP solvers. Applying this strategy to both the Efficient Gauss–Newton ALADIN (see Section~\ref{sec: ALADIN}) and the proposed Efficient Adjoint BFGS ALADIN (Algorithm~\ref{alg:ALADIN}) yields their real-time variants: \emph{RT-Efficient Gauss-Newton ALADIN} and \emph{RT-Efficient Adjoint BFGS ALADIN}.

\subsection{Local Convergence Analysis}\label{sec: convergence} 
The convergence analysis of Algorithm~\ref{alg:ALADIN} follows the same line as in \cite[Section~7]{Houska2016}; therefore, only a proof sketch is provided here due to space constraints. From \cite[Theorem~1]{Du2019}, let the initial iterate of Algorithm~\ref{alg:ALADIN} be sufficiently close to the optimal solution $x^*$. Then, under Assumptions~\ref{ass 1} and~\ref{ass 2} and provided that the local NLPs \eqref{ALADIN-step1} are solved exactly, there exists a constant $\kappa < \infty$ such that
\begin{equation}\label{eq: local converge 0}\small
 \left\|\begin{matrix}
     x_i^{[k+1]} - y_i^*\\
      \mu_i^{[k+1]} - \mu_i^*
 \end{matrix}\right\| \overset{\eqref{ALADIN-step1}}{\leq} \kappa 
 \left\|\begin{matrix}
     y_i^{[k]} - y_i^*\\
   \lambda^{[k]} - \lambda^*
 \end{matrix}\right\|.
\end{equation}
Furthermore, using the efficient adjoint BFGS information generated by \eqref{auxiliary vectors updates}--\eqref{eq:local variables}, the following inequality holds:
\begin{equation}\label{eq: local converge 1}\small
	\begin{split}
		&\left\|\begin{matrix}
		y^{k+1}-y^*\\
			\lambda^{[k+1]}-\lambda^*
		\end{matrix}\right\|
        \leq
		\left\|\begin{matrix}
			y^{[k+1]}-y^*\\
			\lambda^{[k+1]}-\lambda^*\\
			\tilde{\mu}^{[k+1]} -\mu^{*}
		\end{matrix}\right\|    \overset{\eqref{auxiliary vectors updates}-\eqref{eq:local variables}}{\leq}     \gamma\left\|\begin{matrix}
			x^{[k+1]}-y^*\\
			\lambda^{[k]}-\lambda^*\\
			 \mu^{[k+1]}-\mu^*
		\end{matrix}\right\|.
	\end{split}
\end{equation}
Here $\gamma$ quantifies the discrepancy between the Hessian approximation generated by \eqref{eq:BFGS hessian update} and the exact Hessian; this error can be made sufficiently small, as guaranteed by \cite{diehl2010adjoint}. Combining \eqref{eq: local converge 0} and \eqref{eq: local converge 1}, we obtain that if $\gamma(\kappa+1) < 1$, the following contraction holds:
\begin{equation}\small
	\begin{split}
		&\left\|\begin{matrix}
		y^{[k+1]}-y^*\\
			\lambda^{[k+1]}-\lambda^*
		\end{matrix}\right\|
      \overset{\eqref{eq: local converge 0},\eqref{eq: local converge 1}}{\leq}  \gamma(\kappa+1)\left\|\begin{matrix}
		y^{[k]}-y^*\\
		\lambda^{[k]}-\lambda^*\\
	\end{matrix}\right\|.
	\end{split}
\end{equation}
This establishes the local convergence of Algorithm~\ref{alg:ALADIN}.
Note that, the above analysis also applies to the algorithm discussed in Section~\ref{sec: Real Time Adjoint ALADIN}.

        \section{Numerical Experiment}
\label{sec: exp}
This section evaluates the proposed \emph{Efficient Adjoint-BFGS ALADIN} and \emph{Efficient Gauss–Newton ALADIN} \cite{Wu2025CDC} methods {\ADD{on the differential drive robots problem (which is a practical nonconvex MHE problem) and also presents their real-time implementations.}}
\subsection{Time-splitting Based MHE Model}\label{subsection: distributed model}
The corresponding system dynamics and observation model are defined by
\begin{equation*}\label{eq:robot model dynamics}\small
f(\hspace{-0.4mm}x_n,\hspace{-0.4mm}u_n\hspace{-0.4mm})\hspace{-1mm}=\hspace{-1.4mm}
\begin{bmatrix}
{\phi_n} \\
{\psi_n}\\
{\theta_n}
\end{bmatrix}\hspace{-0.9mm}+\hspace{-0.3mm}T\hspace{-0.8mm}
\begin{bmatrix}
v_n \cos \theta_n \\
v_n\sin \theta_n \\
\omega_n
\end{bmatrix}\hspace{-0.8mm}, y_n\hspace{-1mm}=\hspace{-1.4mm}\begin{bmatrix}
r \\
\alpha
\end{bmatrix}
\hspace{-1.2mm}=\hspace{-1.4mm}
\begin{bmatrix}
\sqrt{\phi_n^2\hspace{-0.8mm}+\hspace{-0.8mm}\psi_n^2} \\
\arctan\hspace{-0.8mm}\left(\hspace{-1mm}\frac{\psi_n}{\phi_n}\hspace{-1mm}\right)
\end{bmatrix}\hspace{-0.8mm}+\hspace{-0.8mm}
\begin{bmatrix}
\nu_r \\
\nu_\alpha
\end{bmatrix}\hspace{-0.8mm}.
\end{equation*}

Here, the vector $x=(\phi, \psi, \theta)^\top$ involves three state variables, representing the lateral position $\phi$, longitudinal position $\psi$, and orientation angle $\theta$. The observation vector $ y=(r,\alpha)^\top$ consists of the relative range $r$ and bearing $\alpha$, while the control inputs $u = (v, \omega)^\top$ include the linear and angular velocities, $v$ and $\omega$, respectively. The measurement noise terms $\nu_r$ and $\nu_\alpha$ follow Gaussian distributions $\nu_{r} \sim \mathcal{N}(0, \sigma_{r}^{2})$ and $\nu_{\alpha} \sim \mathcal{N}(0,\sigma_{\alpha}^{2})$. 

Based on the above control system, the resulting MHE problem at each time step $l$ is formulated by penalizing the system distribution and the fitting error over a prediction horizon of length $L$, and solved in a receding-horizon manner. The corresponding simplified objective function can be expressed as follows:
\begin{equation}\label{eq: MHE}
\min_{x}\frac{1}{2} \left\lVert x_{l-L} - \hat{x}_{l-L} \right\rVert_{P^{-1}}^2\hspace{-0.1cm}+\hspace{-0.1cm}\frac{1}{2} \sum_{n = l-L}^l \left\|h(x_n) - y_n \right\|_{V^{-1}}^2 
\end{equation}
where $\hat{x}_{l-L}$ represents the prior state estimate, $P\in\mathbb{R}^{|x_{n}|\times |x_{n}|}$ and $V\in\mathbb{R}^{|y_{n}|\times |y_{n}|}$ are the covariance matrix of the initial state estimation error and measurement noise, respectively, and $h:\mathbb{R}^{|x_n|}\to \mathbb{R}^{|y_n|}$ denotes the nonlinear measurement function.

To apply the proposed algorithm, we adopt the time-splitting strategy from \cite{Wu2025CDC}), which partitions the prediction horizon into $N$ consecutive sub-windows. For each sub-window indexed by $i$, auxiliary boundary variables 
$z_i = \left((z_i^a)^\top, (z_i^b)^\top\right )^\top$ are introduced to represent the initial and terminal states, respectively. Specifically, $z_i^a$ denotes the initial state $z_i^a = x_{l-L+(i-1)t}$, while a new auxiliary variable $z_i^b$ corresponds to the terminal state coupled with the next sub-window via the constraint $z_i^b = z_{i+1}^a$. Collectively, the optimization variables of the $i$-th sub-problem are defined as $X_i =\left ((z_i^a)^\top, (\tilde{x}_{(i)})^\top, (z_i^b)^\top\right )^\top, X_i \in \mathbb{R}^{|X_i|}$, and the corresponding local objective is given as, for $i = 1$ and $i=N$,
\begin{equation}\label{eq:sub-MHE1}\small
\left  \{\begin{split}
&J_1(X_1)\hspace{-1mm} = \hspace{-1mm}  \frac{1}{2} \| z_1^a \hspace{-1.2mm}-\hspace{-0.7mm} \hat{x}_{l-L} \|_{P^{-1}}^2  
\hspace{-1mm} + \hspace{-1mm}\frac{1}{2}\hspace{-1mm}\sum_{j=l-L}^{l-L+t-1}\hspace{-1.2mm} \|  h(x_j)\hspace{-1mm} - \hspace{-0.7mm}y_j\|_{V^{-1}}^2   ,\\
&J_N(X_N) = \frac{1}{2} \sum_{j = l-L+(N-1)t}^{l} \left\|h(x_j) - y_j \right\|^2_{V^{-1}},
\end{split}\right.
\end{equation}
for $i = 2, \cdots, N-1$,
\begin{equation}\label{eq:sub-MHE2}\small
  \begin{split}
J_i(X_i) = \frac{1}{2} \sum_{j=l-L+(i-1)t}^{l-L+it-1} \|  h(x_j) - y_j\|_{V^{-1}}^2.
\end{split}
\end{equation}
Similarly, the nonlinear dynamic equality constraints are partitioned into independent sub-vectors according to the same time-splitting scheme, for $i = 1, \cdots, N-1$,
\begin{equation}
 \small
\mathcal{F}_i(X_i)\hspace{-1mm} =\hspace{-1.5mm}
\begin{bmatrix}
x_{l-L+(i-1)t+1}\hspace{-0.5mm} -\hspace{-0.5mm} f(z_i^a, u_{l-L+(i-1)t}) \\
x_{l-L+(i-1)t+2}\hspace{-0.5mm} -\hspace{-0.5mm} f(x_{l-L+(i-1)t+1}, u_{l-L+(i-1)t+1}) \hspace{-0.8mm}\\
\vdots \\
z_i^b\hspace{-0.5mm} -\hspace{-0.5mm} f(x_{l-L+it-1}, u_{l-L+it-1})
\end{bmatrix},
\end{equation}
for $i = N$:
\begin{equation}\label{eq: N-dynalic}
\small
\mathcal{F}_i(X_i)\hspace{-1mm} =\hspace{-1.5mm}
\begin{bmatrix}
x_{l-L+(i-1)t+1}\hspace{-0.5mm} -\hspace{-0.5mm} f(z_i^a, u_{l-L+(i-1)t}) \\
x_{l-L+(i-1)t+2} \hspace{-0.5mm}-\hspace{-0.5mm} f(x_{l-L+(i-1)t+1}, u_{l-L+(i-1)t+1}) \\
\vdots \\
x_l\hspace{-0.5mm} - \hspace{-0.5mm}f(x_{l-1}, u_{l-1})
\end{bmatrix}.
\end{equation}

Consequently, based on \eqref{eq:sub-MHE1}-\eqref{eq: N-dynalic}, the MHE problem formulated under the time-splitting scheme can be represented as a general distributed optimization problem, as shown in \eqref{eq: DIS OPT}, where the consistency constraints $z_i^b = z_{i+1}^a$ are equivalently expressed as affine coupling constraints in \eqref{eq: DIS OPT}. The explicit form of $A_i$ can be found in \cite{Wu2025CDC}).

The prediction horizon is set to $25$, divided into $4$ parallel sub-windows. The initial position and orientation of the robot are initialized as $x_0=(0.1,0.1,0.0)^\top$, and the reference trajectory $x^*=(\phi^*,\psi^*,\theta^*)^\top$ is generated by MPC under the same control model. The MHE decision variables are initialized at $(\phi^*,\psi^*,0)^\top$. In the proposed \emph{Efficient Adjoint BFGS ALADIN} setup, the penalty parameter is set to $\rho=25$, and the dual variables $\lambda$ and $\mu$ are initialized to zero. At the first iteration, each node sets $H_i=I$ and initializes $C_i$ with the exact Jacobian.
All simulations are executed on a desktop equipped with a $2.5$ GHz Intel $\text{i}5\text{-}13490\text{F}$ processor ($10$ cores, $16$ threads) and $16$ GB of RAM, using \texttt{Casadi-3.6.6} with \texttt{IPOPT} in \texttt{MATLAB R2024b} on Windows $11$.

\subsection{Numerical Results}

\begin{figure}[htbp]
\centering
\includegraphics[width=0.52\textwidth,height=0.25\textheight]{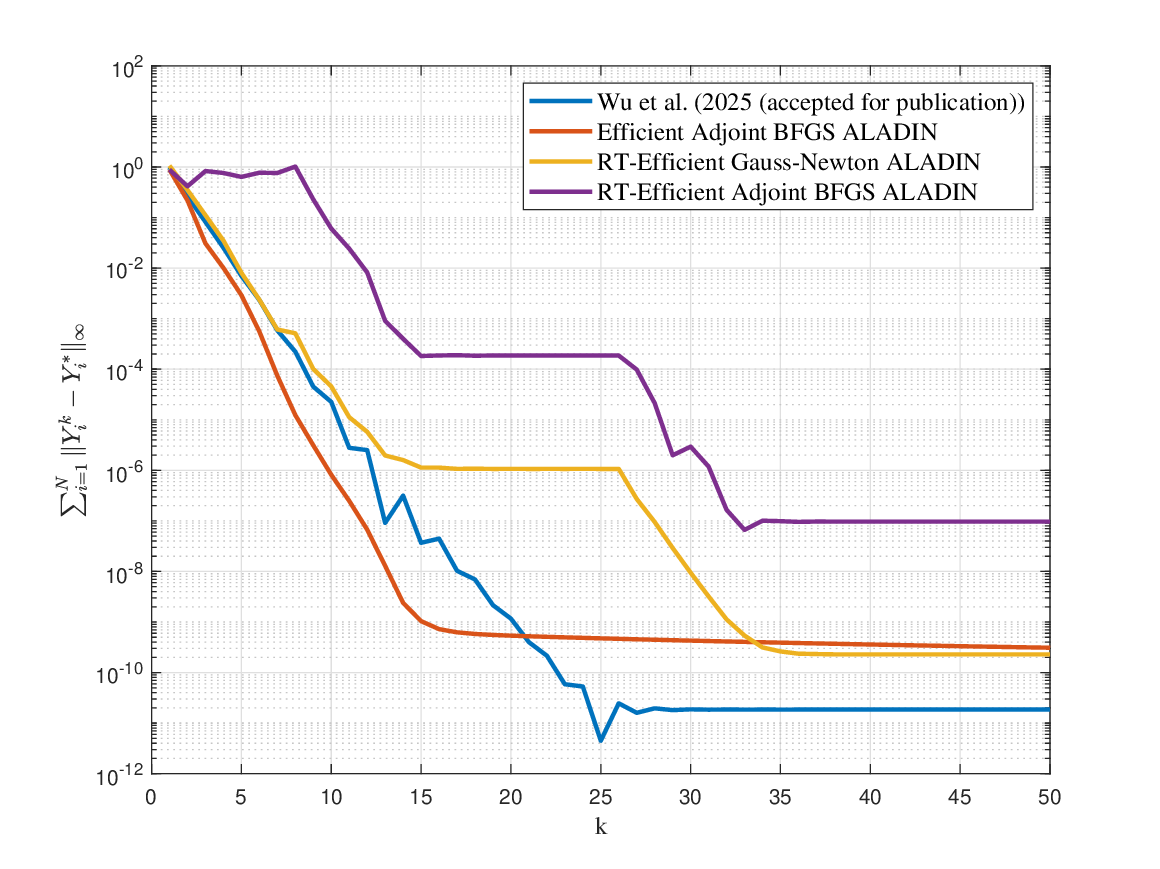}
\caption{Convergence comparison of Efficient Adjoint BFGS ALADIN (Algorithm~\ref{alg:ALADIN}) with Efficient Gauss–Newton ALADIN (in \cite{Wu2025CDC}), RT-Efficient
Gauss-Newton ALADIN and RT-Efficient Adjoint BFGS
ALADIN (in Section~\ref{sec: Real Time Adjoint ALADIN}).} 
\label{fig: convergence plot}
\end{figure}

{\ADD{Fig.}} \ref{fig: convergence plot} compares the convergence behavior of the proposed \emph{Efficient Adjoint BFGS ALADIN}, \emph{Efficient Gauss–Newton ALADIN}, and their real-time variants for solving the MHE problem with $N=4$. Both real-time implementations {\ADD{advertise}} a $\text{log}_3 (k)\in \mathbb N$-based sensitivity update strategy, corresponding to three updates within $50$ iterations. All four methods exhibit linear convergence during sensitivity updates and achieve high-accuracy solutions. Efficient Adjoint–BFGS attains $10^{-9}$ accuracy in $15$ iterations, while \emph{Efficient Gauss–Newton ALADIN} reaches $10^{-11}$ in $25$ iterations. Their real-time counterparts, with fewer sensitivity updates, converge in a two-stage pattern to $10^{-9}$ and $10^{-7}$ within approximately $35$ iterations.

Table~\ref{tab:CPU time} provides the total CPU time of the two efficient ALADIN variants and their real-time counterparts under different sub-window partitions. Since the sub-problem solutions in ALADIN rely on external solvers, their runtimes are nondeterministic and negligible for large-scale cases. Hence, only the CPU time for QP construction and solution in the coordination step is reported. The results show that real-time implementations significantly reduce computation time. Among all methods, \emph{Efficient Gauss–Newton ALADIN} exhibits the longest runtime, while its real-time counterpart achieves the shortest. Unlike the asymptotic relation $N^{*}\approx\sqrt{L}$ reported in (\cite{Wu2025CDC}), the optimal sub-window partition here deviates from the trend because the real-time variants skip sensitivity updates and matrix inversions at each iteration.

\begin{table}[htbp]
\centering
\small
\renewcommand{\arraystretch}{1}
 \scalebox{0.78}{
\begin{tabular}{ 
  |>{\centering\arraybackslash}m{0.4cm}
  |>{\centering\arraybackslash}m{1.8cm}  
  |>{\centering\arraybackslash}m{1.8cm} 
  |>{\centering\arraybackslash}m{1.8cm}  
  |>{\centering\arraybackslash}m{1.8cm}|
  }
\hline
{\centering\textbf{$N$}\rule{0pt}{1.2em}}
& \textbf{Efficient Gauss-Newton ALADIN} 
& \textbf{Efficient Adjoint BFGS ALADIN} 
& \textbf{RT-Efficient Gauss-Newton ALADIN} 
& \textbf{RT-Efficient Adjoint BFGS ALADIN} \\
\hline
\textbf{3}\rule{0pt}{1.2em} & 3.000 & 1.40 & 0.268 & \textbf{0.486} \\
\hline
\textbf{4}\rule{0pt}{1.2em} & 2.300 & 0.971 & 0.232 & 0.516 \\
\hline
\textbf{5}\rule{0pt}{1.2em} & \textbf{1.400} & 0.973 & 0.205 & 0.633 \\
\hline
\textbf{6}\rule{0pt}{1.2em} & 1.500 & \textbf{0.821} & \textbf{0.199} &  0.727\\
\hline
\end{tabular}
}
\caption{Total CPU time $[ms]$ for different algorithms across $N$ sub-windows (measured over $50$ iterations).}
\label{tab:CPU time}
\end{table}

\section{Conclusions} \label{sec:concl}
This paper {\ADD{proposed}} a novel distributed algorithm, termed \emph{Real-Time Efficient Adjoint BFGS ALADIN}, for resource allocation problems. The method efficiently handles non-convex optimization while ensuring linear-scaling communication overhead. Furthermore, we {\ADD{presented}} an event-triggered mechanism that avoids updating local Hessian and Jacobian matrices at every iteration, thereby reducing computational burden. Numerical simulations {\ADD{validated}}  the algorithm's advantages over existing approaches.
	
		
	
    \bibliography{ifacconf,syscop}             
	
	
\end{document}